\documentclass[12pt]{article}

\usepackage{amsfonts}

\newtheorem{thm}{Theorem}[section]

\newtheorem{cor}[thm]{Corollary}

\newtheorem{prop}[thm]{Proposition}

\newcommand{\proof
}{\par\medskip\noindent {\bf Proof.\ \ }}

\newcommand{\be}{\begin{equation}}
\newcommand{\ee}{\end{equation}}

\newcommand{\cG}{\mbox{$\cal G$}}

\newcommand{\cV}{\mbox{$\cal V$}}

\newcommand{\A}{\mathbb A}

\newcommand{\N}{\mathbb N}
\newcommand{\Z}{\mathbb Z}
\newcommand{\Q}{\mathbb Q}

\newcommand{\F}{\mathbb F}

\title{The number of permutations with $k$ inversions}
\author{G\'abor Heged\"{u}s
\\{\normalsize Johann Radon Institute for Computational and Applied Mathematics}
}

\begin{document}
\footnotetext{Research supported in part by OTKA grant K77476}


\maketitle

\begin{abstract}
Let $n\geq 1$, $0\leq t\leq {n \choose 2}$ be arbitrary integers. Define the numbers $I_n(t)$ as the number of permutations of $[n]$ with $t$ inversions.

Let $n,d\geq 1$ and $0\leq t\leq (d-1)n$ be arbitrary integers. Define
{\em the polynomial coefficients} $H(n,d,t)$ as the numbers of compositions of $t$ with at most $n$ parts, no one of which is greater than $d-1$.

In our article we give explicit formulas for the numbers $I_n(t)$ and $H(n,d,t)$ using the theory of Gr\"obner bases and free resolutions.

\end{abstract}
\medskip

\section{ Introduction}

First we introduce some notation. Throughout the paper $n$ is a positive integer. We denote by $[n]$ the set $\{1,\dots ,n\}$. 

We write ${n \choose k}_0$ for the restricted binomial coefficients
\[
{n \choose k}_0:=\left\{ \begin{array} {ll}
{n \choose k}, & \textrm{if $n\geq k$}; \\
\\
0, & \textrm{otherwise.} 
\end{array} \right. 
\]

Let $n\geq 1$ and $d\geq 1$ be arbitrary integers. Throughout the paper $L_n$ stands for the set of vectors 
$$
L_n:=\{(\alpha_1,\ldots ,\alpha_n)\in {\mathbb N}^n:~  0\leq \alpha_i\leq i-1 \mbox{ for }1\leq i \leq n\},
$$
whereas we write $M(n,d)$ for the following set of vectors: 
$$
M(n,d):=\{(\alpha_1,\ldots ,\alpha_n)\in {\mathbb N}^n:~  0\leq \alpha_i\leq d-1 \mbox{ for }1\leq i \leq n\}.
$$

In the following $S_n$ denotes the {\em complete symmetric group}.

Let $\pi=(a_1,a_2,\ldots ,a_n)$ be a permutation of the set $[n]$. Then the pair $(a_i,a_j)$ is an {\em inversion} of the permutation $\pi$, if $i<j$ and $a_i > a_j$. For example, the permutation $(3,1,4,2)$ has three inversions $(3,1)$, $(3,2)$, $(4,2)$.

Let $n\geq 1$, $0\leq t\leq {n \choose 2}$ be arbitrary integers. Define
$$
I_n(t):=|\{\pi\in S_n:~ \pi \mbox{ has }t\mbox{ inversions}\}|.
$$

We give an explicit formula for the number of permutations of $[n]$ with $t$ inversions $I_n(t)$ and for the number of compositions of $t$ with at most $n$ parts, no one of which is greater than $d-1$, $H(n,d,t)$.  E. Netto and D. Knuth gave explicit formulas for $I_n(t)$ only in the restricted case $t\leq n$, but as we know, it did not appear such general explicit formulas for these numbers in the literature. We use sophisticated methods in our proof as the theory of Gr\"obner basis and free resolutions.

We collect here some important fundamental properties of $I_n(t)$. \\

\begin{thm} (\cite[Section 6.4 Theorem D (IV)]{C}) Let $n\geq 1$, $0\leq t \leq {n \choose 2}$ be arbitrary integers. Then
$$
I_n({n \choose 2}-t)=I_n(t).
$$ 
\end{thm}

\begin{thm} (Muir, 1898, \cite{M2})
Denote by 
\begin{equation} \label{gendef}
G_n(z):=\sum_{t=0}^{{n \choose 2}} I_n(t)z^t\in {\Q}[z] 
\end{equation}
the generating function of the numbers $I_n(t)$. Then
\begin{equation} \label{genfg2}
G_n(z)=\frac{(1-z)\cdot\ldots\cdot (1-z^n)}{(1-z)^n}=\prod_{i=1}^{n-1} \Big(\sum_{j=0}^i z^j\Big)
\end{equation}
for each $n\geq 1$.
\end{thm}

\begin{thm} (Bourget, 1871, \cite[Section 6.4 Theorem B]{C})
The number $I_n(t)$ of permutations of $[n]$ with $t$ inversions satisfies the following recurrance relations: 
$$
I_n(t)=\sum_{max(0,\ t-n+1)\leq j\leq t} I_{n-1}(j)
$$ 
for each $n\geq 1$.\\
$I(n,0)=1$ for each $n\geq 1$ and  $I(0,k)=0$ if $k\geq 1$.\\
\end{thm}

E. Netto \cite[ p. 96 ]{N} and D. Knuth \cite[p. 15]{Kn1} gave the following explicit formula for $I_n(t)$ in the case $t\leq n$:
\begin{thm} (\cite[p. 51 (2.51)]{B})
Let $n\geq 1$, $0\leq t\leq n$ be integers. Then
$$
I_n(t)={n+t-1 \choose t}+\sum_{j=1}^{\infty} (-1)^j{n+t-u_j-j-1 \choose t-u_j-j}+
$$
\begin{equation}
+\sum_{j=1}^{\infty}(-1)^j{n+t-u_j-1 \choose t-u_j},
\end{equation}
where \\
$$
u_j:=\frac{j(3j-1)}{2}
$$
are the pentagonal numbers.
\end{thm}

Let $n,d\geq 1$ and $0\leq t\leq (d-1)n$ be arbitrary integers. Define the
{\em polynomial coefficients} $H(n,d,t)$ as
\begin{equation}
H(n,d,t):=|\{(\alpha_1, \ldots ,\alpha_n)\in {\mathbb N}^n:~ \sum_j \alpha_j =t,\ 0\leq \alpha_i\leq d-1\}|.
\end{equation}
These are the numbers of compositions of $t$ with at most $n$ parts, no one of which is greater than $d-1$. We recall the following well--known properties of the numbers $H(n,d,t)$.

\begin{thm} (\cite[p. 77, Ex. 16]{C}, \cite{A} and \cite[p. 154, Ex. 11]{Ri}) \\
The numbers $H(n,d,t)$ have the following recurrance relations:
$$
H(n,d,t)=H(n-1,d,t)+\ldots +H(n-1,d,t-d+1)
$$ 
for each $n,d\geq 1$, $0\leq t\leq (d-1)n$. \\
$H(n,d,0)=1$ for each $n,d\geq 1$ and $H(1,d,t)=1$ for each $d\geq 1$, $0\leq t\leq d-1$. $H(n,d,t)=0$ if either $t<0$ or $t>n(d-1)$. \\

\end{thm}
First L. Euler \cite{E} described the generating function of the polynomial coefficients. \\
\smallskip
\begin{thm} (\cite[Chapter 1, Ex. 16]{C})
The numbers $H(n,d,t)$ have the following generating function:
$$
(1+x+\ldots +x^{d-1})^n=\sum_{t=0}^{n(d-1)} H(n,d,t)x^t.
$$
\end{thm}

Finally we introduce the notion of $q$-binomial coefficients.

The {\em $q$-binomial coefficient} $\Big[ { n \atop m} \Big]_q$ is a $q$-analog for the binomial coefficient, also called a Gaussian coefficient or a Gaussian polynomial. The  $q$-binomial coefficient is given by 
\begin{equation} \label{Gaussco}
 \Big[ { n \atop m}\Big]_q:=\frac{[n]_q!}{[n-m]_q!\cdot[m]_q!}
\end{equation}
for $n,m \in \N$, where $[n]_q!$ is the $q$-factorial (\cite{K}, p. 26)
$$
[n]_q!:=(1+q)\cdot (1+q+q^2)\cdots (1+q+q^2+\ldots +q^{n-1}).
$$
Clearly we have $\Big[ { n \atop k} \Big]_q=\Big[ { n \atop n-k} \Big]_q$. If we substitute $q=1$ into (\ref{Gaussco}), then this substitution reduces this definition to that of binomial coefficients. 

The organisation of this paper is the following. In Section 2 we collected some preliminaries concerning Gr\"obner bases in polynomial rings and free resolutions of graded modules. In Section 3 we stated our main formulas about the numbers $I_n(t)$ and $H(n,d,t)$. In Section 4 we proved these formulas. Finally in Section 5 we gave a lattice-theoretic interpretation of these results.

\section{Preliminaries}
\subsection{Gr\"obner basis and Hilbert function}

We recall now some basic facts
concerning Gr\"obner bases and standard monomials in polynomial rings.

Let $S$ and $R$ stand for the polynomial rings $\Q[x_1,\ldots ,x_n]$ and $\Q[x_0,\ldots ,x_n]$, respectively. 

We denote by $S_{\leq t}$ the vector space of polynomials $f\in S$ of total degree at most $t$, together with $0$. 

We say that a total order $\prec$ on the monomials of $S$ is a {\em term order}, if 1 is the
minimal element of $\prec$, and $uw\prec vw$ holds for any monomials
$u,v,w$ with $u\prec v$. 
Two important term orders are the lexicographic
order $\prec_l$ and the deglex order $\prec _{dl}$.
The lex order is similar to the order of words in the dictionary. Let $u=x_1^{i_1}x_2^{i_2}\cdots x_n^{i_n}$ and
$v=x_1^{j_1}x_2^{j_2}\cdots x_n^{j_n}$ be two monomials. Then $u$ is smaller
than $v$ with respect to lex iff $i_k<j_k$ holds for the smallest index $k$ such
that $i_k\not=j_k$. Similarly, $u$ is smaller
than $v$ with respect to deglex iff either
$\deg (u)< \deg (v)$, or $\deg(u)=\deg(v)$ and $u\prec_l v$. We have $x_n\prec x_{n-1}\prec\ldots \prec
x_1$, for both lex and deglex.

A {\em degree--compatible} term order orders first by total degree: if $\mbox{deg}(u)<\mbox{deg}(v)$, then $u \prec v$. For example, the deglex order is a degree--compatible term order, whereas lex is not.

The {\em leading monomial} ${\rm lm}(f)$
of a nonzero polynomial $f\in S$ is the largest
(with respect to $\prec$) monomial
which appears with nonzero coefficient in $f$ when written as a linear
combination of different monomials.

Let $I$ be an ideal of $S$. A finite subset $G\subseteq I$ is a {\it
Gr\"obner basis} of $I$ if for every $f\in I$ there exists a $g\in G$ such
that ${\rm lm}(g)$ divides ${\rm lm}(f)$. Since $\prec $ is a well founded
order, it can be shown that $G$ is
actually a basis of $I$, i.e.,  $G$ generates $I$ as an ideal of $S$. 
The following Theorem is very fundamental.

\begin{thm} (\cite[Chapter 1, Corollary
3.12]{CCS} or \cite[Corollary 1.6.5, Theorem 1.9.1]{AL}) Every
nonzero ideal $I$ of $S$ has a Gr\"obner basis.
\end{thm}

A Gr\"obner basis $\{ g_1,\ldots ,g_m\}$ of $I$ is {\em reduced}
if the coefficient of ${\rm lm}(g_i)$ is 1, and no nonzero monomial
in $g_i$ is divisible by any ${\rm lm}(g_j)$, $j\not= i$.
By a theorem of Buchberger (\cite[Theorem 1.8.7]{AL}) a nonzero ideal
has a unique reduced Gr\"obner basis.

A monomial $w\in S$ is called a {\it standard monomial for $I$} if
it is not a leading monomial of any $f\in I$. Let ${\rm Sm}(I,\prec)$ stand
for the
set of all standard monomials of $I$ with respect to the term-order
$\prec$ over $\Q$.
The definition and existence of
Gr\"obner bases (\cite[Chapter 1, Section 4]{CCS}) implies that the set
${\rm Sm}(I,\prec)$ is a basis of the $\Q$-vector-space $S/I$ for all
nonzero ideal $I$. 

Let $\cV\subseteq \Q^n$ be a finite set of points. To study the polynomial functions on $\cV$,
it is natural to consider the ideal $I(\cV)$ corresponding to $\cV$:
$$ 
I(\cV):=\{f\in S:~f(v)=0 \mbox{ whenever } v\in \cV\}. 
$$

Recall that a polynomial $h\in R$ is {\em homogeneous} of total degree $d$, if every term appearing in $h$ has total degree exactly $d$. We denote by $R_t$ the vector space of homogeneous polynomials $f\in R$ of degree $t$, together with $0$.

Let $g(x_1,\ldots ,x_n)\in S$ be a polynomial of total degree $d$. We can expand $g$ as the sum of its homogeneous components $g=\sum_{i=0}^d g_i$, where $g_i$ has total degree $i$. Then
$$
g^h(x_0,\ldots ,x_n):=\sum_{i=0}^d g_i(x_1,\ldots ,x_n)x_0^{d-i}\in R
$$
is a homogeneous polynomial of total degree $d$ in $R$. We say that $g^h$ is the {\em homogenization} of $g$.

An ideal $J\unlhd R$ is said to be {\em homogeneous}, if for each $h\in J$, the homogeneous components $h_i$ of $h$ are in $J$ as well. It can be shown easily (\cite[Chapter 8 , Theorem 3.2]{CLO}) that $J$ is a homogeneous ideal iff $J$ is generated by homogeneous polynomials. 

Let $I\unlhd S$ be an arbitrary ideal in $S$. We define the {\em homogenization of $I$} to be the ideal 
$$
I^h:=\langle f^h:~ f\in I\rangle\subseteq R,
$$ 
where $f^h$ denotes the homogenization of the polynomial $f$. It is easy to verify that the homogenization $I^h$ is a homogeneous ideal in $R$ for any ideal $I\unlhd S$.

\begin{thm} \label{homogen} (\cite[Chapter 8, Theorem 4.4]{CLO})
Let $I\unlhd S$ be an ideal of the polynomial ring $S$. Suppose that $\cG=\{g_1,\ldots ,g_k\}$ is a Gr\"obner basis of $I$  with respect to a degree--compatible term order $\prec$. Then the homogenization of $\cG$, $\cG^{h}:=\{g_1^h,\ldots ,g_k^h \}$ is a Gr\"obner basis of the homogeneous ideal $I^h \unlhd R$. 
\end{thm}

Let  $I\unlhd S$ be an arbitrary ideal of $S$. We define the {\em affine Hilbert function} of $S/I$ as
$$
h_{S/I}^a(t)=\mbox{dim}_{\Q}\ S_{\leq t} / (S_{\leq t}\cap I).
$$

The following Theorem connects the affine Hilbert function to certain sets of standard monomials.

\begin{thm} (Macaulay, \cite[Chapter 9, Proposition 3.4]{CLO}) \label{stdmon}
Let $\prec$ be an arbitrary degree--compatible term order on the monomials of $S$. Then
$$
h_{S/I}^a(t)=|\{x^{\alpha}\in \mbox{Sm}(I,\prec):~ \mbox{deg}(x^{\alpha})\leq t \}|.
$$
\end{thm}

Let  $J\unlhd R$ be a homogeneous ideal. The {\em projective Hilbert function} $h_{R/J}^p(t)$ of the quotient ring $R/J$ is
$$
h_{R/J}^p(t)=\mbox{dim}_{\Q}\ R_t /(R_t\cap J).
$$
We can discover a strong connection between the affine Hilbert function of an ideal $I$ of $S$ and the projective Hilbert function of the homogenization $I^h$.
 
\begin{thm} \label{Hilfg} (\cite[Chapter 9 Theorem 3.12]{CLO})
Let $I\unlhd S$ be an ideal of $S$. Then $h_{R/I^{h}}^p(t)=h_{S/I}^a(t)$ for all $t\geq 0$.
\end{thm}

Finally we use the following terminology. If $x^{\alpha}\in R$ (respectively $S$) is an arbitrary monomial and $F\in R$ (respectively $S$) is an arbitrary polynomial, then we denote by $[x^{\alpha}]F$ the coefficient of $x^{\alpha}$ in $F$.

\subsection{Free resolutions of graded modules}

We collect here some basic facts from the theory of free resolutions.

First we give a natural grading structure on the ring $R=\Q[x_0,\ldots ,x_n]$. This comes from the direct sum decomposition
\begin{equation}
R=\bigoplus_{s\geq 0} R_s 
\end{equation}
into the additive subgroups $R_s=\Q[x_0,\ldots ,x_n]_s$, consisting of the homogeneous polynomials of total degree $s$, together with $0$.

We introduce some terminology for describing free resolutions.
 
We say that $M$ over $R$ is a {\em graded module} with a family of subgroups 
 $\{M_t:~ t\in {\Z}\}$ of the additive group, where 
$M_t$ are the homogeneous elements of degree $t$, if we can write $M$ in the form
$$
M=\bigoplus_{t\in \Z} M_t
$$
and 
$$
R_sM_t \subseteq M_{s+t}
$$
for all $s\geq 0$ and $t\in \Z$.
If $M$ is finitely generated, then $M_t$ are finite dimensional vector spaces over $\Q$. 

\begin{prop}
Let $M$ be a graded $R$-module and let $d \in \Z$ be an arbitrary integer. Let 
$$
M(d):=\bigoplus_{t \in \Z} M(d)_t,
$$
where $M(d)_t:=M_{d+t}$. Then $M(d)$ is again a 
 graded $R$-module. 
\end{prop}
We can consider graded free modules of the form  $R(d_1)\oplus \ldots \oplus R(d_n)$ for any integers $d_1, \ldots , d_n$. We call these free modules the {\em twisted graded free modules}.

Let $M,N$ be graded modules over $R$. Then a $\phi:M \longrightarrow N$ homomorphism 
is a {\em graded homomorphism} of degree $d$,
 if $\phi(M_t) \subseteq N_{t+d}$ for all  $t\in \Z$.
If $M$ is a graded $R$--module, then a {\em graded resolution} of $M$ is a 
resolution of the form  
\begin{equation}
0\longrightarrow F_n \longrightarrow \ldots \longrightarrow F_1 \longrightarrow M \longrightarrow 0,
\end{equation}
where each $F_l$ is a twisted graded free module  
and each homomorphism $\phi_l:F_l \longrightarrow F_{l-1}$ 
is a graded homomorphism of degree zero. \\
\smallskip
The following Theorem is fundamental in the theory of free resolutions. 
 
\begin{thm}
{\bf Graded Hilbert Syzygy Theorem} (\cite[Chapter 6, Theorem 3.8]{CLO2}) Every finitely generated $R$--module has a finite graded resolution of length at most $n$. 
\end{thm}
\smallskip

Suppose that  
\begin{equation} \label{gradd}
0\longrightarrow F_n \longrightarrow \ldots \longrightarrow F_1 \longrightarrow M \longrightarrow 0
\end{equation}
is a graded resolution. We say that the resolution (\ref{gradd}) is {\em minimal} iff 
$\phi_l:F_l \longrightarrow F_{l-1}$ takes the standard basis of 
$F_l$ to a minimal generating set of $im(\phi_l)$ 
for each $l\geq 1$. It can be shown that any two minimal graded resolutions of $M$ are isomorphic.

If $M$ is a finitely generated graded $R$--module,  
then we define the Hilbert function $H_M(t)$ by 
$$
H_M(t):=\mbox{dim}_{\F}\ M_t.
$$

Let $I\unlhd R$ be a homogeneous ideal of $R$. Then it is easy to show that $I$ has a natural graded module structure, set $I_t:=I\cap R_t$. Similarly the quotient ring $R/I$ has a natural grading $(R/I)_t:=R_t/I_t$. Thus it comes out from the definitions that if $M:=R/I$ is the quotient graded $R$-module, then $H_M(t)=h_{R/I}^p(t)$ for each $t\geq 0$.

We can compute easily the Hilbert function of the module $R(d)$.
\smallskip
\begin{prop} Let $d\in \Z$ be an arbitrary integer. Then
$H_{R(d)}(t)={t+d+n \choose n}$ for each $t\in \Z$.
\end{prop}
\smallskip
The following Proposition connects the computation of the Hilbert function to the dimension of the free graded modules in a graded resolution of $M$. 

\begin{prop} \label{alter} (\cite[Chapter 6, Theorem 4.4]{CLO2}) 
Let $M$ be a graded $R$--module. Then for any graded free resolution of $M$  
$$
0\longrightarrow F_n \longrightarrow \ldots \longrightarrow F_1 \longrightarrow M \longrightarrow 0
$$
we have 
$$
H_M(t)=\sum_{j=1}^k (-1)^j \mbox{dim}_{\Q}\ (F_j)_t.
$$
\end{prop}
We can easily compute the Hilbert function of the twisted graded free modules $\bigoplus_{i=1}^{m} R(d_{i})$.

\begin{prop} \label{twist}
Let $F$ be the twisted graded free module $\bigoplus_{i=1}^{m} R(d_{i})$. Then 
\begin{equation}
H_{F}(t)=\sum_{i=1}^{m} {n+d_{i}+t \choose n}.
\end{equation}
\end{prop}

Proposition \ref{alter} and Proposition \ref{twist} implies the following Theorem.
\begin{thm} \label{Hilbert} (\cite[Chapter 6, Proposition 4.7]{CLO2}) Let $M$ be a graded $R$-module with the graded free resolution   
\begin{equation}
0\longrightarrow F_n \longrightarrow \ldots \longrightarrow F_1 \longrightarrow M \longrightarrow 0.
\end{equation}
If each $F_j$ is the twisted free graded module  $F_j=\bigoplus_{i=1}^{m_j} R(d_{i,j})$, then 
\begin{equation}
H_M(t)=\sum_{j=1}^k (-1)^j \sum_{i=1}^{m_j} {n-d_{i,j}+t \choose n}.
\end{equation}
\end{thm}

We call the numbers $m_j$ the {\em Betti numbers} of the module $M$. 

Let $\cV\subseteq {\mathbb P}^n$, $\mbox{dim}\ \cV=k$ be a $k$-dimensional variety in the projective space ${\mathbb P}^n$. The variety $\cV$ is called to a {\em complete intersection}, if its homogeneous ideal $I(\cV)\subseteq R$ 
is generated by $n-k$ elements of $R$. 

If $\cV\subseteq {\mathbb P}^n$ is a complete intersection, then we understand very clearly the minimal graded resolution of the ideal $I(\cV)$. We call the minimal graded resolution of complete intersections the {\em Koszul complex}.

\begin{thm} \label{Koszul} (see \cite[Example 13.16]{JH})
Let $\cV\subseteq {\mathbb P}^n$, $\mbox{dim}\ \cV=n-k$ be a complete intersection of hypersurfaces 
defined by the polynomials $G_1, \ldots ,G_k$, where $G_i$ has degree $d_i$. 

Let $M:=I(\cV)$ denote the graded module of the ideal $I(\cV)$. Then the terms in the minimal graded free resolution of $M$ are
$$
F_i=\bigoplus_{1\leq a_1<\ldots <a_i\leq k} R(-d_{a_1}-\ldots -d_{a_i})
$$
for each $1\leq i\leq n$. The Betti numbers $m_i$ are the binomial coefficients ${n \choose i}$. The homomorphism $\phi_i:F_i \longrightarrow F_{i-1}$ sends the generator $e_{a_1,\ldots ,a_n}$ of the module $F_i$ to the sum $\sum_j (-1)^j G_{a_j}\cdot e_{a_1,\ldots , \widehat{a_j}, \ldots, a_i}$.
\end{thm}

\section{The main results}

In Theorem \ref{main1} and Theorem \ref{main2} we give explicit formulas for the sum of the numbers $H(n,d,t)$ and $I_n(t)$.

\begin{thm} \label{main1}
Let $n,d\geq 1$, $0\leq t\leq n(d-1)$ be arbitrary integers. Then 

\begin{equation}
\sum_{i=0}^t H(n,d,i)=\sum_{i=0}^n (-1)^i {n \choose i}{t-di+n \choose n}_0.
\end{equation}

\end{thm}
\begin{thm} \label{main2}
Let $n\geq 1$, $0\leq t \leq {n \choose 2}$ be arbitrary integers. Then

\begin{equation}
\sum_{i=0}^t I_n(i)={n+t \choose t}+\sum_{j=1}^n \sum_{k={j+1 \choose 2}}^{n+1 \choose 2} (-1)^j  [q^{k-{j+1 \choose 2}}]\Big[ { n \atop j} \Big]_q\cdot {t-k+n \choose n}_0,
\end{equation}
where
$$
[q^{k-{j+1 \choose 2}}]\Big[ { n \atop j} \Big]_q
$$
denotes the coefficient of $q^{k-{j+1 \choose 2}}$ in the Gaussian polynomial $\Big[ { n \atop j} \Big]_q$.
\end{thm}

The following two Corollaries are our main results.

\begin{cor} \label{maincor1}
Let $n\geq 1$, $d\geq 1$ and $0\leq t\leq n(d-1)$ be integers. Then 
\begin{equation}
H(n,d,t)=\sum_{i=0}^n (-1)^i {n \choose i}\Big({t-di+n \choose n}_0-{t-di+n-1 \choose n}_0\Big).
\end{equation}
\end{cor} 
\begin{cor} \label{maincor2}
Let $n\geq 1$, $0\leq t\leq {n\choose 2}$ be integers. 
For $1\leq j\leq k$ and ${j+1 \choose 2}\leq k\leq {n+1 \choose 2}$ we denote by 
$$
C(n,j,k):=[q^{k-{j+1 \choose 2}}]\Big[ { n \atop j} \Big]_q
$$
the coefficient of $q^{k-{j+1 \choose 2}}$ in the  Gaussian polynomial 
$\Big[ { n \atop j} \Big]_q$. 
Then 
$$
I_n(t)={n+t \choose t}-{n+t-1 \choose t-1}+
$$
$$
\sum_{j=1}^n \sum_{k={j+1 \choose 2}}^{n+1 \choose 2} (-1)^j  C(n,j,k) \Big({t-k+n \choose n}_0-{t-k+n-1 \choose n}_0\Big).
$$
\end{cor}

\section{Proofs}

{\bf Proof of Theorem \ref{main1}:} \\
Denote by $\cV$ the {\em $d$-box} $\cV:= [d]^n\subseteq {\Q}^n$ . We can describe easily the reduced Gr\"obner basis and the standard monomials of the ideal $I(\cV)$.
\begin{thm} \label{doboz}
Let $\prec$ be an arbitrary term order 
on the monomials of $\Q[x_1,\dots,x_n]$.
Then the reduced Gr\"obner basis of the ideal $I(\cV)$ is
$$
\cG=\{\prod_{j=1}^d (x_i-j):~ 1\leq i\leq n\}
$$
and the set of standard monomials is 
$$
\mbox{Sm}(I(\cV),\prec)=\{x_1^{\alpha_1}\cdot\ldots\cdot x_n^{\alpha_n}:~ 0\leq \alpha_i\leq d-1\}.
$$
\end{thm}
Hence we infer from Theorem \ref{stdmon} that
$$
h_{S/I(\cV)}^a(m)=|\{x^{\alpha}\in \mbox{Sm}(I(\cV),\prec_{dl}):~ \mbox{deg}(x^{\alpha})\leq m \}|=
$$
\begin{equation} \label{vektor}
=|\{(\alpha_1,\ldots ,\alpha_n)\in M(n,d):~ \sum_i \alpha_i\leq m\}|.
\end{equation}
It is easy to verify that the ideal $I(\cV)$ is a complete intersection in the affin space ${\A}^n$. Then we infer from Theorem \ref{homogen} that the homogeneous ideal $I(\cV)^h$ is a complete intersection in the projective space ${{\mathbb P}}^n$, too.

Thus we can describe the minimal graded free resolution 
\begin{equation} \label{graddd}
0\longrightarrow F_n \longrightarrow \ldots \longrightarrow F_1 \longrightarrow M \longrightarrow 0,
\end{equation}
where $M:=I(\cV)^h$ is the graded module of the homogenization $I(\cV)^h$. The resolution (\ref{graddd}) is the Koszul complex by Theorem  \ref{Koszul}. Here $d_i=d$ for each $1\leq i\leq n$.

We infer from Theorem \ref{Hilfg} that
$$
h_{R/I^h}^p(t)=h_{S/I}^a(t)
$$
for each $t\geq 0$. Hence using the exact sequence 
$$
0\longrightarrow I \longrightarrow R \longrightarrow R/I \longrightarrow 0
$$
we get that
$$
h_{R/I^h}^p(t)=h_R^p(t)-h_{I^h}^p(t).
$$
Finally Theorem \ref{Hilbert} and Theorem \ref{Koszul} implies that 
\begin{equation} \label{vegso}
h_{S/I}^a(t)=h_{R/I^h}^p(t)=h_R^p(t)-h_{I^h}^p(t)=\sum_{i=0}^n (-1)^i {n \choose i}{t-di+n \choose n}_0.
\end{equation}
Theorem \ref{maincor1} follows from the equations (\ref{vektor}) and (\ref{vegso}). 

\noindent
{\bf Proof of Theorem \ref{main2}:} \\
Let 
\begin{equation} 
\cV:=\{\pi(1,\ldots ,n) :~ \pi\in S_n\}\subseteq {\Q}^n
\end{equation}
denote the orbit in $\Q^n$ of the vector $(1,\ldots ,n)$ under the permutation action of the complete symmetric group $S_n$.

We prove Theorem \ref{main2} in the following three steps. \\
1. First we give a combinatorial proof of the equality
\begin{equation} 
\sum_{i=0}^t I_n(i)=|\{(\alpha_1,\ldots ,\alpha_n)\in L_n:~ \sum_i \alpha_i\leq t\}|.
\end{equation}
2. Next we prove using Gr\"obner basis theory  that 
\begin{equation}
h_{S/I(\cV)}^a(t)=|\{(\alpha_1,\ldots ,\alpha_n)\in L_n:~ \sum_i \alpha_i\leq t \}|.
\end{equation}
3. Finally we prove using the theory of free resolutions that
\begin{equation}
 h_{S/I(\cV)}^a(t)={n+t \choose t}+\sum_{j=1}^n \sum_{k={j+1 \choose 2}}^{n+1 \choose 2} (-1)^j  [q^{k-{j+1 \choose 2}}]\Big[ { n \atop j} \Big]_q\cdot {t-k+n \choose n}_0,
\end{equation}
where
$$
[q^{k-{j+1 \choose 2}}]\Big[ { n \atop j} \Big]_q
$$
denotes the coefficient of $q^{k-{j+1 \choose 2}}$ in the  Gaussian polynomial $\Big[ { n \atop j} \Big]_q$. \\

1.  It follows from  (\ref{genfg2}) that $G_n(z)$, the generating function of $I_n(t)$, has the following expression:

\begin{equation} \label{genfg3}
G_n(z)=\prod_{i=1}^{n-1} \Big(\sum_{j=0}^i z^j\Big).
\end{equation}

Consider the following polynomial $P$ in $n$ variables:
\begin{equation} \label{genfg4}
P(x_1,\ldots ,x_n):=\prod_{i=2}^{n} \Big(\sum_{j=0}^{i-1} x_i^j\Big)\in \Z[x_1, \ldots ,x_n].
\end{equation}
Then the equations (\ref{genfg3}) and (\ref{genfg4}) yield to the equation  
\begin{equation} \label{pgenfg}
P(z,\ldots ,z)=G_n(z). 
\end{equation}
The definition of the polynomial $P$ implies that
\begin{equation} \label{genfg5}
P(x_1,\ldots ,x_n)=\sum_{\alpha\in L_n} x^{\alpha}, 
\end{equation}
where if $\alpha$ is the vector $(\alpha_1,\ldots , \alpha_n)$, then we use the notation $x^{\alpha}:=x_1^{\alpha_1}\cdots x_n^{\alpha_n}$.

This implies that 
$$
|\{(\alpha_1,\ldots ,\alpha_n)\in L_n:~ \sum_i \alpha_i\leq t \}|=
$$
$$
=\sum_{j=1}^t |\{x^{\alpha}:~ \alpha\in L_n,\ \mbox{deg}(x^{\alpha})=j \}|=
$$
\begin{equation} \label{Psum}
 =\sum_{j=1}^t |\{x^{\alpha}:~ [x^{\alpha}]P=1,\  \mbox{deg}(x^{\alpha})=j\}|,
\end{equation}
where we used the expansion (\ref{genfg5}) in the equality (\ref{Psum}). Clearly we have
$$
\sum_{j=1}^t |\{x^{\alpha}:~ [x^{\alpha}]P=1,\   \mbox{deg}(x^{\alpha})=j\}|=\sum_{j=1}^t [z^j]P(z,\ldots ,z),
$$
where $[z^j]P(z,\ldots ,z)$ denotes the coefficient of $z^j$ in the polynomial $P(z,\ldots ,z)$.  Finally we infer from (\ref{gendef}) and (\ref{pgenfg}) that
$$
\sum_{j=1}^t [z^j]P(z,\ldots ,z)=\sum_{j=1}^t [z^j] G_n(z)=\sum_{j=1}^t I_n(j).
$$

2. We recall the definition of the complete and elementary symmetric polynomials.
Let $i$ be a nonnegative integer. Write 
$$ 
h_i(x_1,\ldots,x_n)=\sum_{a_1+\cdots +a_n=i}x_1^{a_1}x_2^{a_2}\cdots
x_n^{a_n}
$$
for the $i$-th complete symmetric polynomial.
Clearly $h_i\in {\Q}[x_1,\ldots ,x_n]$ is the sum of all monomials of total degree
$i$. 

For $0\leq i\leq n$ we write $\sigma_i$ for the $i$-th elementary 
symmetric polynomial:
$$ \sigma_i(x_1,\ldots,x_n)=\sum_{S\subset [n],~|S|=i} x_S.$$ 
$\sigma_i\in \Q[x_1,\ldots ,x_n]$ is the sum of all squarefree monomials 
of degree $i$ in the variables
$x_1,\ldots,x_n$.

For $1\leq k\leq n$ we introduce the polynomials $f_k\in S$ as follows:
$$
f_k=\sum_{i=0}^k (-1)^i h_{k-i}(x_k,x_{k+1},\ldots,x_n)\sigma_i
(1,\ldots,n).
$$

Clearly $f_k\in \Q[x_k,x_{k+1},\ldots, x_n]$. Moreover, 
$\mbox{deg}(f_k)=k$ and the 
leading monomial of $f_k$ is $x_k^k$ with respect to any term order $\prec$
for which $x_1\succ x_2\succ \ldots \succ x_n$. 

G. Heged\"us and L. R\'onyai \cite[Theorem 2.2]{HNR} described the reduced Gr\"obner basis and the standard monomials of $I(\cV)$ with respect to any term order 
on the monomials of ${\Q}[x_1,\dots,x_n]$.
\begin{thm} 
\label{perm}
Let $\prec$ be an arbitrary term order 
on the monomials of $\Q[x_1,\dots,x_n]$ such that $x_n\prec \ldots \prec x_1$.
Then the reduced Gr\"obner basis of $I(\cV)$ is
$$
\{f_i:~ 1\leq i \leq n\}.
$$ 
Moreover the set of standard monomials is 
\begin{equation}
\label{standard1}
\mbox{\em Sm}(I(\cV),\prec)=\{x_1^{\alpha_1}\ldots x_n^{\alpha_n}:~ 
0\leq \alpha_i \leq i-1, ~\mbox{\em for}~1\leq i\leq n \}.
\end{equation}

\end{thm}

Hence we infer from Theorem \ref{stdmon} that
\begin{equation}
h_{S/I(\cV)}^a(t)=|\{x^{\alpha}\in \mbox{Sm}(I(\cV),\prec_{dl}):~ \mbox{deg}(x^{\alpha})\leq t \}|=
\end{equation}

\begin{equation}
=|\{(\alpha_1,\ldots ,\alpha_n)\in L_n:~ \sum_i \alpha_i\leq t\}|
\end{equation}
for each $t\geq 0$. \\
3. Theorem \ref{perm} implies  that the zero--dimensional ideal $I(\cV)$ is a complete intersection in the affin space ${\A}^n$. But then Theorem \ref{homogen} gives that the homogenization $I(\cV)^h$ is a complete intersection in the projective space ${{\mathbb P}}^n$. 
Let $M:=I(\cV)^h$ denote the graded module of the homogenization  $I(\cV)^h$. Then 
the minimal graded resolution 
\begin{equation}
0\longrightarrow F_n \longrightarrow \ldots \longrightarrow F_1 \longrightarrow M \longrightarrow 0,
\end{equation}
is the Koszul complex, where $d_i=i$ for each $1\leq i \leq n$.

Let $n\geq 1$, $1\leq j\leq n$ and $1\leq k\leq {n+1 \choose 2}$ be arbitrary integers. Define the numbers
\begin{equation}
D(n,j,k):=|\{I\subseteq [n]:~ |I|=j,\ \sum_{i\in I} i=k\}|.
\end{equation}
Clearly $D(n,j,k)=0$ if $1\leq k< {j+1 \choose 2}$ for a fixed $1<j\leq n$.

Theorem \ref{Hilbert} and Theorem \ref{Koszul} implies that
$$
h_{S/I}^a(t)=h_{R/I^h}^p(t)=h_R^p(t)-h_{I^h}^p(t)=
$$
\begin{equation} \label{Hilbert2}
={n+t \choose t}+\sum_{j=1}^n \sum_{k={j+1 \choose 2}}^{n+1 \choose 2} (-1)^j  D(n,j,k)\cdot {t-k+n \choose n}_0.
\end{equation}

We give in the following a simpler expression for the numbers $D(n,j,k)$.
\begin{thm} \label{Gauss} (\cite[Chapter 2, Theorem 2.25]{B})
Let $n\geq 1$, $1\leq j\leq n$ and $1\leq k \leq {n+1 \choose 2}$. Then
\begin{equation} \label{Dgauss}
D(n,j,k)=[q^{k-{j+1 \choose 2}}]\Big[ { n \atop j} \Big]_q,
\end{equation}
where $[q^{k-{j+1 \choose 2}}]\Big[ { n \atop j} \Big]_q$ denotes the coefficient of $q^{k-{j+1 \choose 2}}$ in the  Gaussian polynomial $\Big[ { n \atop j} \Big]_q$.
\end{thm}

We include here a proof of Theorem \ref{Gauss} for the reader's convenience.

{\bf Proof of Theorem \ref{Gauss}:}

We can recognise the numbers $D(n,j,k)$ as the coefficients in certain two--variables polynomials.
\begin{thm} \label{coeff}
Let $n\geq 1$ be an arbitrary integer. 
Define
$$
Q_n(t,q):=\prod_{i=1}^n (1+tq^i)\in {\Z}[t,q].
$$
Then 
$$
D(n,j,k)=[q^kt^j]Q_n(t,q),
$$
the coefficient of $[q^kt^j]$ in the polynomial $Q_n(t,q)$, for each $1\leq j\leq k$ and $1\leq k\leq {n+1 \choose 2}$.
\end{thm}
\proof

This is clear from the definition of the polynomial $Q_n(t,q)$ and the numbers $D(n,j,k)$.

The following well--known Theorem is a special case of the $q$-binomial Theorem (\cite[Theorem 1.1]{G} and \cite[Chapter 2, Exercise 11]{C}). 

\begin{thm} ({\bf Cauchy binomial Theorem}) \label{Cauchy}
Let $n\geq 1$ be an arbitrary integer. Then 
\begin{equation}
Q_n(t,q)=\prod_{i=1}^n (1+tq^i)=\sum_{m=0}^n t^m q^{\frac{m(m+1)}{2}}\Big[ { n \atop m} \Big]_q.
\end{equation}
\end{thm}

Finally we conclude from Theorem \ref{coeff} and Theorem \ref{Cauchy} that
$$
D(n,j,k)=[q^kt^j]Q_n(t,q)=[q^kt^j]\Big(\sum_{m=0}^n t^m q^{\frac{m(m+1)}{2}}\Big[ { n \atop m} \Big]_q\Big)=
$$
$$
=[q^kt^j]\Big(t^jq^{\frac{j(j+1)}{2}}\Big[ { n \atop j} \Big]_q\Big),
$$
which implies (\ref{Dgauss}). $\Box$

Now we finish the proof of Theorem \ref{main2}. From equations (\ref{Hilbert2}) and (\ref{Dgauss}) we infer that
$$
h_{S/I}^a(t)={n+t \choose t}+\sum_{j=1}^n \sum_{k={j+1 \choose 2}}^{n+1 \choose 2} (-1)^j  D(n,j,k)\cdot {t-k+n \choose n}_0=
$$
\begin{equation}
={n+t \choose t}+\sum_{j=1}^n \sum_{k={j+1 \choose 2}}^{n+1 \choose 2} (-1)^j  [q^{k-{j+1 \choose 2}}]\Big[ { n \atop j} \Big]_q\cdot {t-k+n \choose n}_0.
\end{equation}

\section{Concluding remarks}

We give a nice interpretation of this problem and its solution in lattice theory. 

We say that a $P$ poset is {\em graded of rank $n$} if every maximal chain of $P$ has the same length. In this case we can give a unique rank function $\rho:P \to \{0, \ldots ,n\}$ such that $\rho(x)=0$ if $x$ is a minimal element of $P$, and $\rho(y)=\rho(x)+1$ if $y$ covers $x$ in $P$. If $\rho(x)=i$, then we say that $x$ has {\em rank} $i$. 

If $P$ is a graded poset of rank $n$ and has $t_i$ elements of rank $i$ then the polynomial 
$$
G(P,q)=\sum_{i=0}^n t_iq^i\in \Z[q]
$$ 
is called the {\em rank-generating function} of $P$.

Let $n>0$ be an integer. We make the set of all positive divisors of $n$ easily into a poset $D_n$ by defining $i\preceq j$ in $D_n$ if $j$ is divisible by $i$. It is easy to verify that $D_n$ with this partial ordering $\preceq$ becomes a graded poset. The rank of $P$ is the number of prime divisors of $n$, counting with multiplicity.

Let $P_n$ denote the partially ordered set of all $n$-permutations in which $p\leq q$ if either $p$ can be obtained from $q$ by a series of transpositions, or $p=q$. Then $\leq$ is called the Bruhat order on $S_n$ and it is easy to verify that $(P_n,\leq)$ is again a graded poset.

Let $p_1, \ldots ,p_{n}$ denote the first $n$ prime numbers in increasing order. Then it can be shown easily that the poset $(P_n,\leq)$ is isomorphic with the poset $(D_k,\preceq)$, where $k=p_1^1 \cdot \ldots \cdot p_{n-1}^{n-1}$. Thus $I_n(t)$ is the number of elements of rank $t$ in $D_k$. Similarly, if $\ell:=p_1^{d-1}\cdot \ldots \cdot p_n^{d-1}$, then $H(n,d,t)$ is the number of elements of rank $t$ in $D_{\ell}$. 

Thus we determined the rank generating function of the ranked poset $D_n$ for certain $n$'s. It would be very interesting to determine the rank generating function of $D_n$ in the general case, when $n$ is an arbitrary integer. 

R. S. Deodhar and M. Srinivasan \cite{DS} defined a statistic on involutions, which they called {\em weight}.

Let $F(2n)$ denote the set of all fixed point free involutions of $[2n]$. 

A $2$-cycle is a set consisting of two distinct positive integers. We write $[i,j]$ for an arc with $i<j$. The {\em span} of an arc $[i,j]$ is defined as $\mbox{span}[i,j]=j-i-1$. A pair of disjoint arcs $[i,j]$ and $[k,l]$ is a {\em crossing} if $i<j<k<l$ or $k<i<l<j$.

We can write the involutions in {\em standard representations}, which is in increasing order in their initial points.
 
Let $\delta$ be an arbitrary involution. The numbers of arcs in $\delta$ is denoted by $|\delta|$. The {\em crossing number} of $\delta$, denoted by $c(\delta)$, is the number of pairs of arcs of $\delta$ that are crossings. Define the {\em weight} of $\delta$, denoted by $\mbox{wt}(\delta)$, as follows:
$$
\mbox{wt}(\delta):=\Big(\sum_{[i,j]\in \delta}   \mbox{span}[i,j]\Big)-c(\delta).
$$

Let $\delta=[a_1,b_1][a_2,b_2]\ldots [a_n,b_n]\in F(2n)$ be an arbitrary involution. Then we say that $\tau$ is obtained from $\delta$ by an {\em interchange}, written $\delta \sim \tau$, if there exist $1\leq i<j \leq n$ such that \\
\\
(i)\ $\tau$'s standard representation is obtained from $\delta$ by exchanging $b_i$ and $a_j$, or\\
(ii)\ $\tau$'s standard representation is obtained from $\delta$ by exchanging $b_i$ and $b_j$.

We say that $\tau$ is obtained from $\delta$ by a {\em weight increasing interchange}, if $\delta \sim \tau$ and $\mbox{wt}(\delta)<\mbox{wt}(\tau)$.

They defined a partial order on $F(2n)$ as follows: Let $\delta,\tau\in F(2n)$. Then $\delta \leq \tau$ if $\tau$ can be obtained from $\delta$ by a sequence of (zero or more) weight increasing interchanges.

R. S. Deodhar and M. Srinivasan \cite[Theorem 1.3]{DS} proved  that $(F(2n), {\leq})$ is a graded poset of rank $2{n\choose 2}$ for each $n\geq 1$. The rank of $\delta \in F(2n)$ is given by $\mbox{wt}(\delta)$ and the rank generating function of this poset $F(2n)$ is the following:
$$
G(F(2n),q)=[1]_q[3]_q\ldots [2n-1]_q.
$$
We think that it would be interesting to apply our methods for this poset and we can give explicit formulas for the rank generating function of $F(2n)$, but we shall treat this topic in a future paper.

\medskip
\noindent
{\bf Acknowledgements.} I am indebted to Jonathan Farley for his useful remarks.

\end{document}